\documentclass[reqno,ceqn]{amsart}
\usepackage{amssymb,amscd,verbatim}
\usepackage{enumerate}
\numberwithin{equation}{section}
\usepackage{bm}
\usepackage{color}
\usepackage{mathrsfs}
\usepackage{bigints}
\usepackage{cases}
\usepackage{cancel}
\usepackage{ulem}
\usepackage{enumerate}
\usepackage{graphicx}
\usepackage{url}
\usepackage{nccmath}

\newtheorem{thm}{Theorem}[section]
\newtheorem{cor}[thm]{Corollary}
\newtheorem{lem}[thm]{Lemma}

\theoremstyle{definition}

\newtheorem{rem}{Remark}[section]

    {\begin{thm}\label{#1} \ifShowLabels \TeXref{#1} \fi}%
    {\end{thm}}

    {\begin{def}\label{#1} \ifShowLabels \TeXref{#1} \fi}%
    {\end{def}}

    {\begin{lem}\label{#1} \ifShowLabels \TeXref{#1} \fi}%
    {\end{lem}}

    {\begin{cor}\label{#1} \ifShowLabels \TeXref{#1} \fi}%
    {\end{cor}}

\newcommand{\R}{\mathbb{R}}

\begin{document}
	
	\begin{center}
		
		{\large \bf
			Remarks on comparison principles for $p$-Laplacian \\
			with extension to $(p,q)$-Laplacian 
		}
		\vspace{1cm}
		
		{Ahmed Mohammed \\
			Department of Mathematical Sciences \\
			Ball State University \\ Muncie, Indiana 47306, USA\\
			{email: amohammed@bsu.edu}
			\vspace{0.5cm}
			
			Antonio Vitolo\\
			Dipartimento di Matematica\\ Universit\`a di Salerno\\84084 Fisciano (Salerno),
			Italy\\
			{email: vitolo@unisa.it}
		}
		
	\end{center}
	
	\vspace{0.5cm}

	\vspace{0.5cm}
	
	\vspace{0.1cm}

	\vspace{1cm}

		\noindent{\bf Abstract}:  Our purpose is to generalize some recent comparison principles for operators driven by $p$-Laplacian to a  wide class of quasilinear equations including $(p,q)$-Laplacian. It turns out, in particular, that
		 adding a $q$-Laplacian to $p$-Laplacian allows to weaken the assumptions needed on the Hamiltonian of lower order terms. The results are specialized in the case that the Hamiltonian has at most polynomial growth in the gradient with coefficients depending on $x$ and $u$.
		\vspace{0.5cm}

	\section{\bf Introduction and main results}\label{intro}
	
	\setcounter{page}{1}
Comparison principles are powerful and widely used tools to investigate qualitative properties of solutions of second-order elliptic and parabolic equations. Among numerous works in the literature that deal with comparison principles we just mention  \cite{ABM, AS, BM, BNMP, CL, LPR, Lind, MOV, POR, SEP, VER}. 
 In a recent paper \cite{MOV}, the authors studied comparison principles for the $p$-Laplacian of the form
 \begin{equation}\label{pfg}
	-\Delta_pu+H(u,\nabla u)=0
\end{equation}
in a bounded domain $\Omega$ of $\R^n$, targeted to the case of $H(t,\xi)=f(t)+g(t)|\xi |^\sigma$, with $\sigma >0$, under different assumptions on the Hamiltonian $H$, which complement and extend the results of  \cite{LPR}. For a strong comparison principle we refer to \cite{PRR}. It is also worth remarking that the authors use in \cite{MOV} a different approach with respect to their previous papers \cite{MPV,MV20}, devoted to the Harnack inequality and the strong maximum principle for $p$-Laplacian. 

	In this paper we will focus on the principal term of equation (\ref{pfg}), also including a dependence of $x\in \Omega$. Namely, we consider the more general quasilinear equation
\begin{equation}\label{exp}
	-\textup{div}(\bm{a}(x,\nabla u))+H(x,u,\nabla u)=0,
\end{equation}
where $\bm{a}:\Omega\times\mathbb\R^n\to\mathbb R^n$ is a Carath\'eodory function, i.e. measurable in $x \in \Omega$ and continuous in $\xi \in \R^n$, and $H:\Omega\times\mathbb R\times\mathbb R^n\to\mathbb R$ is a continuous function, in order to prove comparison principles of broad application. Such assumptions are intended to hold throughout all the paper.

Exploiting the major generality, we will see that for a significant class of vector fields  $\bm a$, the assumptions on the Hamiltonian $H$ can be weakened, with respect to the previous results for $p$-Laplacian contained in \cite{MOV}. 
\vskip0.2cm
 In particular, aiming to capture the well-known $(p,q)$-Laplacian, we want to consider vector fields 
\begin{equation}\label{pt}
	\bm{a}(x,\xi)= b(x)|\xi|^{p-2}\xi+c(x)|\xi|^{q-2}\xi, \; \; (x,\xi)\in\Omega\times\mathbb R^n,
\end{equation}
with constants $p\ge q>1$, which encompass $p$-Laplacian, when $q=p$ and $b(x)+c(x)\equiv 1$, and  $(p,q)$-Laplacian when $q<p$ and $b(x)=c(x)\equiv 1$. More generally, we are led to  the following ellipticity assumption: 
\vskip0.2cm
\begin{itemize}
	\item[{\bf (a-e)}] \quad 	$(\bm{a}(x,\xi)-\bm{a}(x,\eta))\cdot(\xi-\eta)\ge \gamma\left\{\begin{array}{ll}
		(1+|\xi|+|\eta|)^{p-2}|\xi-\eta|^2,&\text{if}\;\;1<p<2\\[.3cm]
		|\xi-\eta|^r,&\text{if}\;\;p\ge 2,
	\end{array}\right. $ 
\end{itemize}  
\vskip0.2cm
for all $(x,\xi),(x,\eta)\in\Omega\times\mathbb R^n$, where $\gamma$ is  a positive constant and the latter inequality holds for all $ r \in [\max(2,q),p]$. 
\vskip0.2cm
We will see in the next section that in fact both $p$-Laplacian and $(p,q)$-Laplacian are defined through vector fields $\bm a=\bm a(\xi)$ which meet  the ellipticity condition {\bf (a-e)}.  As a further example, see \cite[Lemma 1]{PTO},  condition {\bf(a-e)} holds for $r=p=q$ if   $\bm{a}(x,\xi)$ belongs to $ C^1(\Omega\times\mathbb R^n\setminus\{0\})\cap C^0(\Omega\times \mathbb R^n)$, and satisfies the following condition for some constants $c>0$ and $0\le\kappa\le 1$:
$$\sum_{i,j=1}^n\frac{\partial\bm{a}_j}{\partial \eta_i}(x,\eta)\xi_i\xi_j\ge c(\kappa+|\eta|)^{p-1}|\xi|^2$$
for all $(x,\eta)\in \Omega\times\mathbb R^n\setminus\{0\}$ and all $\xi\in\mathbb R^n$.
\vskip0.2cm
In recent years, the $(p,q)$-Laplacian has appeared in many investigations. Not even close to be complete, we refer for instance to the works \cite{COF, Rad1, LPG, LNZ, Rad0, Rad2, Rad3, Rad4}, and the references therein.

Concerning $H$, we assume that $H:\Omega \times \R\times \R^n\to \R$ is continuous, and it satisfies the following two conditions:
\vskip0.1cm
{\bf{(H-1)}\;\;$H(x,t,\xi)-H(x,s,\xi)\ge \omega(x,t-s) + \tau(x,t-s) |\xi|^\sigma$}
\vskip0.1cm
\noindent for all $(x,\xi)\in \Omega \times\mathbb R^n$ and $s,t\in \R$ such that $s<t$, where $\omega, \tau:\Omega\times(0,\infty)\to \R$ are non-negative,   lower semicontinuous in $x\in \Omega$ and non-decreasing in $r\in (0,\infty)$;
\vskip0.1cm
{\bf{(H-2)} \;\;$|H(x,t,\eta)-H(x,t,\xi)|\le M||\eta|^\sigma-|\xi|^\sigma|$}
\vskip0.1cm
\noindent for all $(x,t) \in K \times \R$ and $\xi,\eta \in \R^n$, where $K$ is any compact subset of $\Omega$, with positive constants $\sigma$ and $M=M(K)$, or in a weaker  version, for all $(x,t) \in K \times J$ and $\xi,\eta \in \R^n$, where $K\times J$ is any compact subset of $\Omega\times \R$, with positive constants $\sigma$ and $M=M(K,J)$

	\vskip0.2cm
	Note that here we are simplifying the assumptions with respect to our previous paper \cite{MOV}, see Section 1, to minimize the complications due to the additional dependence on $x$ and to avoid too many technicalities. Still comparing with \cite{MOV}, we point out that, since the use of exponents $(p,q)$ is standard for double Laplacians, the positive exponent $q$ of \cite{MOV} corresponds here to $\sigma$. A further technical complexity arises from having three exponents instead of two, and so we have made an effort to smoothen the arguments in our proofs. In particular, in the case $q=p$, they provide an alternative, lighter proof of the results in Section 3 of \cite{MOV} for $p$-Laplacian.  We will also state the corresponding results in the case of special Hamltonians $H(x,u,\nabla u)=f(x,u)+g(x,u)|\nabla u|^\sigma$.
	\vskip0.2cm
	Now, let us write
		\begin{equation}\label{qq}
			Qu:=-\textup{div}(\bm{a}(x,\nabla u)) +H(x,u,\nabla u),
		\end{equation}
	
For $u,v\in W^{1,p}_{loc}(\Omega)$, $1 < p <\infty$,  the inequality $Qu\le Qv$ in $\Omega$ means that
		\begin{equation}\label{def1}
			\int_\Omega(\bm{a}(x,\nabla u)-\bm{a}(x,\nabla v))\cdot \nabla\varphi\le \int_\Omega(H(x,v,\nabla v)-H(x,u,\nabla u))\varphi
		\end{equation}
		for all $\varphi\in C^1_c(\Omega)$ such that $\varphi\ge 0\;\text{in}\;\Omega.$ To have the integrals well-defined, we need some consistency assumptions. For instance, in addition to the previous assumptions on $H$, we require
			\begin{equation}\label{intp}
				|\bm{a}(x,\nabla u)|, \,|\bm{a}(x,\nabla v)|, \,H(x,u,0), \,H(x,v,0)  \in L^{\frac p{p-1}}_{loc}(\Omega).
			\end{equation}
		  These  also work  whenever test functions $\varphi \in W^{1,p}(\Omega)$, with compact support in $\Omega$, are used in (\ref{def1}). If $u,v\in W^{1,\infty}_{loc}(\Omega)$, it is sufficient to require
		  	\begin{equation}\label{intinfty}
		  		|\bm{a}(x,\nabla u)|, \,|\bm{a}(x,\nabla v)|, \,H(x,u,0), \,H(x,v,0)  \in L^{1}_{loc}(\Omega).
		  	\end{equation}
	\vskip0.1cm
	For more specific assumptions on $H$, we refer to \cite{MOV}, where it is also shown how  such assumptions imply that the right-hand side in (\ref{def1}) is well defined. A discussion of this aspect is provided in the proof of Lemma 2.1 there, and we have chosen to avoid such digression here. 
	
	We will always intend that consistency assumptions like (\ref{intp}), or (\ref{intinfty}) at the occurrence, are satisfied, whenever we use the inequality $Qu \le Qv$, in order that the integrals in (\ref{def1}) are well defined when using test functions $\varphi \in W^{1,p}(\Omega)$ with compact support in $\Omega$. By $u\le v$ in $\Omega$, we mean that 
	$$u(x) \le v(x) \; \; \text{for a.e. } x \in \Omega.$$ 
	On the other hand, inequality $u\le v$ on $\partial\Omega$ is to be understood in the sense that
	$$\limsup_{x\to\partial\Omega}(u(x)-v(x))\le 0.$$
	
		\vspace{.1cm}
		
Now we establish the first main result in our general setting. It is a generalization of \cite[Theorem 1.1]{MOV} for $p$-Laplacian and includes $(p,q)$-Laplacian. One of the highlights is that the presence of the additional $q$-Laplacian allows to weaken the assumptions on the Hamiltonian $H$ in certain cases.

	\begin{thm}\label{cp:thm}  Assume $2 \le  p<\infty$ and $1<q\le p$. Let conditions {\bf{(a-e)}} and  {\bf{(H-i)}} be satisfied, $i=1,2$. Suppose  $u,v \in W^{1,p}_{loc}(\Omega)$, and one of the following conditions:
		\vskip0.2cm
		\begin{description}
			\item[A$_1$] \; \,$0<\sigma\le 1$, $\omega>0$; \vspace{.1cm}
			\item[A$_{2a}$]  \, $1\le \sigma\le  \max(q,2)-1$, $\omega >0$ and $\tau>0 $; \vspace{.1cm}
			\item[A$_{2b}$]  \;  $ \max(q,2)-1\le \sigma\le p-1$, $\tau>0$;		\vspace{.1cm}
			\item[A$_3$]  \;  \;$p-1\le  \sigma\le p$, $\tau >0$ and $\textup{ess}\sup_{\Omega}(u-v)<\infty$.
		\end{description}
		\vskip0.1cm
		Suppose $u, v\in W^{1,\infty}_{loc}(\Omega)$, and 
		\begin{description}
			\item[A$_4$]  \, \,	$p\le\sigma$, $\tau>0$.
			\end{description}
		\vskip0.1cm
		If $Qu \le Qv$ in $\Omega$, and $u \le v$ on $\partial \Omega$, then  $u\le v$ in $\Omega$.
	\end{thm}
	
	\vskip0.1cm
	\begin{rem}{\rm We point out that, if
	\begin{equation}\label{essup}
	 \textup{ess}\sup_{\Omega}(u-v)< \infty,
	\end{equation}
	as for for instance in the case that $u,v \in L^\infty_{loc}(\Omega)$, it is sufficient to adopt the weaker version of condition {\bf{(H-2)}}. This is the case for assumptions  {\bf A$_3$} and {\bf A$_4$}  of Theorem \ref{cp:thm}. The same is possible for {\bf A$_1$}, {\bf A$_{2a}$} and {\bf A$_{2b}$},  if we add  condition (\ref{essup}). \qed}
	\end{rem}

	\vskip0.2cm
	With respect to the previous results for $p$-Laplacian, see Theorem 1.1 of \cite{MOV},  the range $1\le \sigma \le p-1$, where $\omega>0$ is required in the case of $p$-Laplacian,   splits  into two ranges, $1\le \sigma \le \max(q,2)-1$ and $\max(q,2)-1\le \sigma\le p-1$, in the case of $(p,q)$-Laplacian: in the first range we continue to need $\omega>0$, but in the latter one we can relax this condition to the default assumption $\omega\ge 0$. Being $q < p$ for the $(p,q)$-Laplacian, this is an interval of positive length, unless $p=2$, when such interval shrinks to a point. Indeed, this is the rule in the case of $p$-Laplacian. In fact, taking $q=p \ge 2$, the range $\max(q,2)-1\le \sigma\le p-1$ is always reduced only to a point, $\sigma=p-1$. 
	
	\vskip0.2cm
	In the next result, concerning the case $1 < p < 2$, we assume that $u,v \in W^{1,\infty}_{loc}(\Omega)$ for all $\sigma>0$. It is a generalization of \cite[Theorem 1.2]{MOV}, but this time the ranges established for $p$-Laplacian remain unchanged.

	\begin{thm}\label{cp:thmm}  Assume $1 < p < 2$. Let conditions {\bf{(a-e)}} and  {\bf{(H-i)}} be satisfied, $i=1,2$. Suppose  $u,v \in W^{1,\infty}_{loc}(\Omega)$, and one of the following conditions: 
		\begin{description}
			\item[${\bf{\widetilde A}}_1$] \; $0<\sigma<1$, $\omega >0$;\vspace{.1cm}
			\item[${\bf{\widetilde A}}_2$] \; $1\le \sigma$.
		\end{description}
		If $Qu \le Qv$ in $\Omega$, and $u \le v$ on $\partial \Omega$, then  $u\le v$ in $\Omega$.
	\end{thm}
	
		\vskip0.2cm
	
	It is interesting and useful to establish  the counterpart of the above comparison principles for special Hamiltonians 
	\begin{equation}\label{H*}
		H^*(x,s,\xi)=f(x,s)+g(x,s)|\xi|^\sigma
	\end{equation}
	\vskip0.2cm
In this case we suppose basically that the functions $f$ and $g$ are continuous in $\Omega \times \R$. To have  assumptions playing the role of {\bf{(H-1)}} and  {\bf{(H-2)}} for general Hamiltonians, we will assume that $f(x,t)$ and $g(x,t)$ are non-decreasing in $t\in \R$ for all $x \in \Omega$, which correspond to $\omega \ge 0$ and $\tau \ge 0$, respectively. Likewise, assuming $f(x,t)$ and $g(x,t)$ increasing in $t\in \R$ for all $x \in \Omega$ corresponds to suppose $\omega> 0$ and $\tau > 0$, respectively.  See the discussion in Section \ref{inequalities} with particular reference to (\ref{fgxeq}), (\ref{fpos}) and (\ref{gpos}).

	\vskip0.2cm
	
	In view of this, we can state the following comparison principles for our quasi-linear elliptic operators with special Hamiltonians $H^*$. In this case we write
		\begin{equation}\label{qq*}
			Q^*u:=-\textup{div}(\bm{a}(x,\nabla u)) +f(x,u)+g(x,u)|\nabla u|^\sigma,
		\end{equation}
	where $\bm a $ is a Carathéodory function in $\Omega \times \R^n$, $f$ and $g$ continuous functions in $\Omega \times \R$. 
	
	Here are the results.

	\begin{thm}\label{cpfg:thm}  Assume $2 \le  p<\infty$ and $1<q\le p$. Let condition {\bf{(a-e)}} be satisfied, and $f(x,t)$, $g(x,t)$ be non-decreasing function in $t \in \R$ for all $x \in \Omega$. Suppose  $u,v \in W^{1,p}_{loc}(\Omega)\cap L^\infty_{loc}(\Omega)$, and one of the following conditions:
		\vskip0.2cm
		\begin{description}
			\item[A$_1$] \, \, $0<\sigma\le 1$, $f$ increasing in $t$; \vspace{.1cm}
			\item[A$_{2a}$]  \;  $1\le \sigma\le  \max(q,2)-1$, $f$ and $g$ increasing in $t$; \vspace{.1cm}
			\item[A$_{2b}$]  \;  $ \max(q,2)-1\le \sigma\le p-1$, $g$ increasing in $t$;		\vspace{.1cm}
			\item[A$_3$]  \;  \;\,$p-1\le  \sigma\le p$, $g$ increasing in $t$, and $\textup{ess}\sup_{\Omega}(u-v)<\infty$.
		\end{description}
		\vskip0.1cm
		Suppose $u, v\in W^{1,\infty}_{loc}(\Omega)$, and
		\begin{description}
			\item[A$_4$]  \; \,	$p\le\sigma$, $g$ increasing in $t$.
		\end{description}
		\vskip0.1cm
		If $Qu \le Qv$ in $\Omega$, and $u \le v$ on $\partial \Omega$, then  $u\le v$ in $\Omega$.
	\end{thm}
	
	\vskip0.2cm

	\begin{thm}\label{cpfg:thmm}  Assume $1 < p < 2$. Let condition {\bf{(a-e)}} be satisfied, and $f(x,t)$, $g(x,t)$ be non-decreasing function in $t \in \R$ for all $x \in \Omega$. Suppose  $u,v \in W^{1,\infty}_{loc}(\Omega)$, and one of the following conditions: 
		\begin{description}
			\item[${\bf{\widetilde A}}_1$] \; $0<\sigma<1$, $f$ increasing;\vspace{.1cm}
			\item[${\bf{\widetilde A}}_2$] \; $1\le \sigma$.
		\end{description}
		If $Qu \le Qv$ in $\Omega$, and $u \le v$ on $\partial \Omega$, then  $u\le v$ in $\Omega$.
	\end{thm}
	
	\vskip0.2cm
	The paper is organized as follows. In Section \ref{inequalities} we show that our setting comprises $(p,q)$-Laplacian and discuss the interplay between general and special Hamiltonians. In Section \ref{comparison} we give proof of the comparison principles stated above and discuss variants which can be useful in the applications.

	\section{\bf Motivations and basic inequalities}\label{inequalities}
	
	 We start this section showing that the ellipticity condition {\bf{(a-e)}} holds for $p$-Laplacian 
	 \begin{equation}\label{p-L}
	 	\Delta_pu:= \text{div} (|\nabla u|^{p-2}\nabla u)
	 \end{equation}
	 with $p>1$ and $(p,q)$-Laplacian 
	 	\begin{equation}\label{pq-L}
	 		\Delta_{(p,q)}u:= \Delta_pu+\Delta_qu
	 	\end{equation}
	 with $p>q>1$.
	 \vskip0.1cm
	 While checking this, all the inequalities are assumed to hold for all $\xi,\eta\in \R^n$.
	 \vskip0.1cm
	 In the case of $p$-Laplacian, in which $\bm{a}(x,\xi) = |\xi|^{p-2}\xi$, we know by \cite[Chapter 12]{Lind} that 
	 \begin{equation}\label{aepL}
	 		(|\xi|^{p-2}\xi-|\eta|^{p-2}\eta)\cdot(\xi-\eta) 
	 	\ge \gamma_p\left\{\begin{array}{ll}
	 		(1+|\xi|+|\eta|)^{p-2}|\xi-\eta|^2,&\text{if}\;\;1<p<2\\[.1cm]
	 		|\xi-\eta|^p,&\text{if}\;\;p\ge 2,
	 	\end{array}\right.
	 \end{equation}
	 for all $\xi,\eta\in \R^n$, where $\gamma_p$ is a positive constant, and so condition {\bf{(a-e)}} is satisfied for $r=p=q$. 
	 Note in particular that for all $p>1$
	 \begin{equation}\label{pLpos}
	 	(|\xi|^{p-2}\xi-|\eta|^{p-2}\eta)\cdot(\xi-\eta) \ge 0.
	 \end{equation}
	 Moreover,  for $p\ge 2$, the constant $\gamma_p$ can be taken non-increasing with $p$.

	 In the case of $(p,q)$-Laplacian  it is
	 \begin{equation}\label{pqLxi}
	 	\bm{a}(x,\xi)= |\xi|^{p-2}\xi+|\xi|^{q-2}\xi,
	 \end{equation}
	 and we observe that by (\ref{pLpos}) we have 
	 \begin{equation}\label{pqLpq}
	 (\bm{a}(x,\xi)-\bm{a}(x,\eta))\cdot(\xi-\eta) \ge	(|\xi|^{r-2}\xi-|\eta|^{r-2}\eta)\cdot(\xi-\eta)
	 \end{equation}
	 for both $r=q$ and $r=p$.
	
	 Now, we consider three ranges:  (i) $p>q \ge 2$, \; (ii) $p>2 > q>1$, \; (iii) $2 \ge  p>q>1$. 
	 \vskip0.1cm
	 In the range (i), by   (\ref{pqLpq}) and the second inequality in (\ref{aepL}), we have
	\begin{equation}\label{Lpq(i)}
			(\bm{a}(x,\xi)-\bm{a}(x,\eta))\cdot(\xi-\eta)\ge \gamma_p\, \max\left(|\xi-\eta|^p,|\xi-\eta|^q\right) \ge\gamma_p\, |\xi-\eta|^r
	\end{equation}
for all $r \in [q,p]$. So condition {\bf{(a-e)}} holds with $\gamma=\gamma_p$ in the range (i). The same can be obtained in the range (iii), using  (\ref{pqLpq}) with $r=p$ and the first inequality in (\ref{aepL}).

In order to deal with the range (ii), we go back to the derivation of inequality (\ref{aepL}). In \cite[Chapter 12]{Lind} it is deduced from the inequality
\begin{equation}\label{pre-arepL}
		(|\xi|^{p-2}\xi-|\eta|^{p-2}\eta)\cdot(\xi-\eta) \ge \frac12\left(|\xi|^{p-2}+|\eta|^{p-2}\right)|\xi-\eta|^2,
\end{equation}
which holds for all $p>1$.
From this, we get for the vector field  (\ref{pqLxi}) 
		\begin{equation}\label{pre-aepqL}
			(\bm{a}(x,\xi)-\bm{a}(x,\eta))\cdot(\xi-\eta) \ge \frac12\left(|\xi|^{p-2}+|\xi|^{q-2}+|\eta|^{q-2}+ |\eta|^{p-2}\right)|\xi-\eta|^2.
		\end{equation}

	Now, if  $q < 2 < p$, the function $h(t)=t^{p-2}+t^{q-2} $  is positive, continuous in $(0,\infty)$ and goes to infinity both for $t\to 0^+$ and $t\to \infty$. Then it has a  minimum $m_{p,q}>0$ depending on $p$ and $q$, and so in the case (ii)
		\begin{equation}\label{aepqL2}
			(\bm{a}(x,\xi)-\bm{a}(x,\eta))\cdot(\xi-\eta) \ge m_{p,q} |\xi-\eta|^2.
		\end{equation}
	On the other hand, since $p>2$, by (\ref {pqLpq}) and the second inequality in (\ref{aepL}) we also have 
		\begin{equation}\label{aepqLp}
			(\bm{a}(x,\xi)-\bm{a}(x,\eta))\cdot(\xi-\eta) \ge \gamma_p |\xi-\eta|^p.
			\end{equation}
		So, letting $\gamma=\min(m_{p,q},\gamma_p)$, we have 
	\begin{equation}\label{Lpq(ii)}
		(\bm{a}(x,\xi)-\bm{a}(x,\eta))\cdot(\xi-\eta) \ge \gamma \max(|\xi-\eta|^p,|\xi-\eta|^2) \ge \gamma\,|\xi-\eta|^r
	\end{equation}
 for all $r \in [2,p]$. 
Collecting inequalities (\ref{Lpq(i)}), obtained in the ranges (i) and (iii),  and (\ref{Lpq(ii)}), for the range (ii), we get the full condition {\bf(a-e)} with $\gamma=\gamma(p,q)$.
\vskip0.3cm
Next, consider the Hamiltonian $H:\Omega \times \R \times \R^n \to \R$, which we assume to be continuous. Conditions {\bf{(H-1)}} and {\bf{(H-2)}}  yield:

\begin{equation}\label{Hest}
	\begin{split}
	H(x,s,\xi)-H(x,t,\eta) &=(H(x,s,\xi)-H(x,t,\xi))+(H(x,t,\xi)-H(x,t,\eta))\\
	&\ge \omega(x,s-t)+\tau(x,s-t) |\xi|^\sigma- M||\xi|^\sigma-|\eta|^\sigma|,
	\end{split}
\end{equation}
that is
\begin{equation}\label{Hxeq}
	H(x,t,\eta)-H(x,s,\xi)  +\omega(x,s-t)+\tau(x,s-t) |\xi|^\sigma\le  M||\xi|^\sigma-|\eta|^\sigma|
\end{equation} 
for $s >t$, where by the non-decreasing monotonicity of $\omega(x,r)$ and $\tau(x,r)$ in $r$, we have $\omega(x,s-t) \ge0$ and $\tau(x,s-t)\ge 0$.
\vskip0.1cm
Suppose now that $\omega(x,r) >0$, respectively $\tau(x,r) >0$, in $\Omega \times \R$. If $s-t \ge \lambda_1$, where $\lambda_1$ is a positive constant, and $K$ is a compact subset of $\Omega$, then by the non-decreasing monotonicity in $r$, and the lower semicontinuity in $x$ of $\omega$, respectively $\tau$, we have 
\begin{equation}\label{omegapos}
	\omega(x,s-t) \ge \omega(x,\lambda_1) \ge \omega_1 \; \; \text{for } x \in K \;  \text{and } s-t \ge \lambda_1
	\end{equation}
with a constant $\omega_1>0$, respectively
	\begin{equation}\label{taupos}
		\tau(x,s-t) \ge \tau(x,\lambda_1) \ge \tau_1  \; \; \text{for } x \in K \;  \text{and } s-t \ge \lambda_1
	\end{equation}
with a constant $\tau_1>0$. 
\vskip0.2cm
In the case of special Hamiltonians $H^*(x,s,\xi)=f(x,s)+g(x,s)|\xi|^\sigma$, with $f$ and $g$ continuous functions in $\Omega \times \R$, we see below that inequalities like (\ref{Hxeq}), (\ref{omegapos}) and (\ref{taupos}) hold with $f(x,s)-f(x,t)$ and $g(x,s)-g(x,t)$ in the place of $\omega(x,s-t)$ and $\tau(x,s-t)$, respectively. In this case, the non-decreasing (increasing) monotonicity of $f$ and $g$ in $s$ plays the role of the non-negativity (positivity) of $\omega$ and $\tau$, respectively.

In fact, let $K\subset \Omega$ and $J \subset \R$ be compact. Arguing as above for (\ref{Hxeq}), we get
\begin{equation}\label{fgxeq}
	\begin{split}
		 H^*(x,t,\eta)\hskip-0.05cm -\hskip-0.05cm H^*(x,s,\xi)  \hskip-0.05cm+\hskip-0.05cm f(x,s)\hskip-0.05cm -\hskip-0.05cm f(x,t)\hskip-0.05cm+\hskip-0.05cm(g(x,s)\hskip-0.05cm-\hskip-0.05cm g(x,t)) |\xi|^\sigma\le  M||\xi|^\sigma\hskip-0.05cm-\hskip-0.05cm|\eta|^\sigma|
	\end{split}
\end{equation}
with  $M= \max_{K\times J}|g|$, and, assuming that $f(x,s)$ and $g(x,s)$ are non-decreasing in $s$, we have  $f(x,s)-f(x,t) \ge 0$ and $g(x,s)-g(x,t)\ge 0$ for $s>t$.

Suppose now $f(x,s)$, respectively $g(x,s)$,  increasing in $s$. If $s-t \ge \lambda_1>0$ for $s,t \in J$,  then  by the increasing monotonicity in $s\in \R$ and the continuity of $f$, respectively $g$, in $\Omega \times \R$, there exists $(x_1,t_1)\in K \times J$ such that
	\begin{equation}\label{fpos}
		f(x,s)-f(x,t) > f(x,t+\lambda_1) -f(x,t) \ge \omega_1 \; \; \text{for } x \in K \;  \text{and } s-t \ge \lambda_1
	\end{equation}
where $\omega_1=  f(x_1,t_1+\lambda_1) -f(x_1, t_1)>0$, respectively
	\begin{equation}\label{gpos}
		g(x,s)-g(x,t) > g(x,t+\lambda_1) -g(x,t) \ge \tau_1 \; \; \text{for } x \in K \;  \text{and } s-t \ge \lambda_1
	\end{equation}
where $\tau_1=  g(x_1,t_1+\lambda_1) -g(x_1, t_1)>0$.  

\section{\bf Proof of the main results and variants}\label{comparison}

In order to prove comparison principles, we will go further with the inequality (\ref{Hxeq}) for the variations of $H$ in Section \ref{inequalities}, which allow to exploit conditions ${\bf A}_\alpha$  of Theorem \ref{cp:thm}. Here is the proof of our main result.

\vskip0.2cm

{\bf Proof of Theorem \ref{cp:thm}}. Let $u,v \in W^{1,p}_{loc}\Omega)$ such that $Qu \le Qv$ in $\Omega$,  and $u \le v$ on $\partial \Omega$. In order to show that $u \le v$ in $\Omega$, we argue by contradiction, supposing on the contrary that
\begin{equation}\label{lamda}
	\lambda:= \textup{ess}\sup_{\Omega}(u-v)>0.
\end{equation}

Then we find an increasing sequence of positive numbers $\lambda_k< \lambda$, $k \in \mathbb N$, such that $\lambda_k\to \lambda$ as $k \to \infty$.  Let $w_k=(u-v-\lambda_k)^+ $. It is a non-negative function in $W^{1,p}(\Omega)$. By virtue of the boundary condition $u \le v$ on $\partial \Omega$, the non-negative functions $w_k$ belong to $W^{1,p}_0(\Omega)$. In addition, letting $\mathcal O_k$ be the set where $u-v > \lambda_k$, we have $\mathcal O_k \subset \mathcal O_1 $, and dist$(\mathcal O_1, \partial \Omega)>0$. So there exists a compact subset $K \subset \Omega$ such that supp\,$w_k \subset K$ for all $k \in \mathbb N$.

By the consistency conditions, see (\ref{intp}) and (\ref{intinfty}), we can use the functions $w_k$ to test the inequality $Qu \le Qv$ in $\Omega$.  Hence
	\begin{equation}\label{t1}
		\int_\Omega(\bm{a}(x,\nabla u)-\bm{a}(x,\nabla v))\cdot \nabla w_k \le \int_\Omega(H(x,v,\nabla v)-H(x,u,\nabla u))w_k.
	\end{equation}
Since $p \ge 2$, from condition {\bf{(a-e)}} we have 
	\begin{equation}\label{t1s}
	\gamma_0	\int_\Omega|\nabla w_k|^r \le\int_\Omega(H(x,v,\nabla v)-H(x,u,\nabla u))w_k
	\end{equation}
for all $r \in [\max(q,2),p]$ and, by  (\ref{Hxeq}), 
	\begin{equation}\label{t1d}
		 \gamma_0\int_\Omega|\nabla w_k|^r  + \int_\Omega \omega(x,u-v)\,w_k+\int_\Omega \tau(x,u-v)\,|\nabla u|^\sigma w_k \le M\int_\Omega||\nabla u|^\sigma-|\nabla v|^\sigma |w_k.
	\end{equation}

Since $w_k=0$ outside a compact $K$ and $u-v \ge \lambda_1$ where $w_k>0$, and $\omega$ is non-decreasing and lower semicontinuous, we have
	\begin{equation}\label{inf}
		\int_\Omega \omega(x,u-v)\,w_k \ge  \int_K \omega(x,\lambda_1)\,w_k \ge \omega_1\int_\Omega w_k,
	\end{equation}
where $\omega_1\ge 0$ is the minimum of $\omega(x,\lambda_1)$ in $K$. The same holds for $\tau$, and we denote by $\tau_1\ge 0$ the minimum of $\tau(x,\lambda_1)$ in $K$. Thus from (\ref{t1d}) we get
	\begin{equation}\label{step0}
	\gamma_0	\int_\Omega|\nabla w_k|^r + \omega_1\int_\Omega w_k+\tau_1\int_\Omega |\nabla u|^\sigma w_k \le M\int_\Omega||\nabla u|^\sigma-|\nabla v|^\sigma |w_k.
	\end{equation}
	
	The next step is to majorize the latter integral in order to obtain the inequality 
	\begin{equation}\label{step1}
		\gamma_0\int_\Omega|\nabla w_k|^r + \omega_1\int_\Omega w_k \le M_1\int_{\Omega_k}|\nabla w_k|^\sigma w_k,
	\end{equation}
where $\Omega_k=\{x \in \mathcal O_k: \nabla u \neq \nabla v\}$, for some positive constant $M_1$.
	
	If $0 < \sigma \le 1$ we use the inequality
	\begin{equation}\label{sig1:eq}
		||\xi|^\sigma - |\eta|^\sigma| \le |\xi-\eta|^\sigma,
	\end{equation}
	which directly  yields (\ref{step1}), in the case ${\bf A}_{1}$, with $M_1=M$. We point out that in this case we do not need condition $\tau>0$, which instead has to be used in the remaining case.
	
	For $\sigma >1$, we use the inequality
\begin{equation}\label{1sig:eq}
	||\xi|^\sigma - |\eta|^\sigma| \le \left(\tfrac{1+\varepsilon}{\varepsilon}\right)^{\sigma-1}\,|\xi-\eta|^\sigma+\frac12\,((1+\varepsilon)^{\sigma-1}-1)\,|\xi|^\sigma,
\end{equation}
which holds for all $\varepsilon>0$ and reduces to (\ref{sig1:eq}) for $\sigma=1$ (see \cite[Lemma A.2]{MOV} for a proof). 
In this case, we exploit condition $\tau>0$, which is assumed in the cases ${\bf A}_{2a}$, ${\bf A}_{2b}$, ${\bf A}_{3}$ and ${\bf A}_{4}$. Estimating  the right-hand side of (\ref{step0}) by (\ref{1sig:eq}) with $\xi=\nabla u$ and $\eta=\nabla v$, we get
	\begin{equation}\label{tauest}
		\begin{split} 
	\gamma_0	\int_\Omega|\nabla w_k|^r + \omega_1\int_\Omega w_k+\tau_1\int_\Omega |\nabla u|^\sigma w_k &\le M\left(\tfrac{1+\varepsilon}{\varepsilon}\right)^{\sigma-1}\int_\Omega|Dw_k|^\sigma w_k\\
	&+ \frac M2\,((1+\varepsilon)^{\sigma-1}-1)\int_\Omega |\nabla u|^\sigma w_k 
\end{split}
\end{equation}
Taking $\varepsilon >0$ small enough in order to have $\frac M2\,((1+\varepsilon)^{\sigma-1}-1)< \tau_1$ and subtracting the second integral of the right-hand side, we  obtain (\ref{step1}) with $M_1=M\left(\tfrac{1+\varepsilon}{\varepsilon}\right)^{\sigma-1}$.

\vskip0.2cm

{\it Case  $0< \sigma < p$ }
\vskip0.1cm

In this case we estimate  the right-hand side of (\ref{step1}) by Young inequality, introducing a parameter $r>\sigma$ to be chosen in the sequel:
	\begin{equation}\label{young}
		\int_\Omega|\nabla w_k|^\sigma w_k = \int_\Omega(\delta|\nabla w_k|^\sigma) (\delta^{-1}w_k) \le \tfrac {\sigma}r \,\delta^r\int_\Omega |\nabla w_k|^r+\tfrac{r-\sigma}{r}\,\delta ^{-\frac r{r-\sigma}} \int_\Omega w_k^{\frac r{r-\sigma}}.
	\end{equation}
Choosing  $\delta>0$ small enough, in order to have $M_1\,\tfrac r{\sigma} \,\delta^r <\gamma_0$, and subtracting from (\ref{step1}) the first integral of the above right-hand side, multiplied by $M_1$,  we get
	\begin{equation}\label{step2}
		\gamma_1\int_\Omega|\nabla w_k|^r + \omega_1\int_\Omega w_k \le M_2\int_{\Omega_k}w_k^{\frac r{r-\sigma}}
	\end{equation}
with positive constants  $\gamma_1=\gamma_0-M_1\tfrac {\sigma}r \,\delta^r$ and $M_2= M_1\tfrac{r-\sigma}{r}\,\delta ^{-\frac r{r-\sigma}} $.

We will show for all cases, from ${\bf A}_1$ to ${\bf A}_3$ (except $\sigma=p$), an inequality of type
\begin{equation}\label{invPS}
	\int_\Omega|\nabla w_k|^r \le C	\int_{\Omega_k}|w_k|^r 
\end{equation}
for some $r \in [\max(q,2),p]$, with a constant $C>0$ independent of $k$.  Postponing its proof, we observe that (\ref{invPS}) leads to a contradiction. In fact, combined with Sobolev inequalities, 
this yields
\begin{equation}\label{ctrd}
\int_{\Omega}|w_k|^{\overline r}\le C\left(\int_\Omega|\nabla w_k|^r \right)^{\frac{\overline r}{r}}\le C	\left(\int_{\Omega_k} w_k^r \right)^{\frac{\overline r}{r}} \le C	|\Omega_k|^{\frac{\overline r}{r}-1}\int_{\Omega}|w_k|^{\overline r}
\end{equation} 
with $\overline r>r$ (the Sobolev conjugate of $r$, if $r< n$, any number greater than $r$, otherwise), where $C$ is used to denote possibly different constants.
\vskip0.1cm
Here the first (and equally the last) integral cannot be zero, otherwise $w_k=0$ a.e. in $\Omega$, so that $u-v \le \lambda_k< \lambda$ a.e. in $\Omega$, against (\ref{lamda}). Since $|\Omega_k|\to 0$ as $k \to \infty$, we get a contradiction, $1 \le C |\Omega_k|^{\frac{\overline r}{r}-1}$, for large $k\in \mathbb N$. We conclude that  $\lambda \le 0$, i.e. \hskip-0.1cm$u \le v$ in $\Omega$, as it was to show.
\vskip0.1cm
We are left with proving (\ref{invPS}), and we will do it range by range, using (\ref{step2}).
\vskip0.1cm

In the cases ${\bf A}_1$ ($0 < \sigma \le 1$, $\omega>0$) and ${\bf A}_{2a}$ ($1\le \sigma\le  \max(q,2)-1$, $\omega>0$ and $\tau>0$), we set $r=p$. Since $\frac p{p-\sigma}\le p$, we can estimate the integral on the right-hand side of  (\ref{step2}) as follows, introducing $\varepsilon>0$ to be chosen in the sequel:
\begin{equation}\label{sig<1}
	\begin{split}
	\int_{\Omega_k}w_k^{\frac{p}{p-\sigma}}
		&\le \int\limits_{w_k < \varepsilon}w_k^{\frac{p}{p-\sigma}}+\int\limits_{w_k \ge \varepsilon}w_k^{\frac{p}{p-\sigma}} \\
		&=\int\limits_{w_k < \varepsilon} w_k^{\frac{\sigma}{p-\sigma}}w_k+\varepsilon^{\frac{p}{p-\sigma}}\int\limits_{w_k \ge \varepsilon} \left(\frac{w_k}{\varepsilon}\right)^{\frac{p}{p-q}}\\
		&\le \varepsilon^{\frac{\sigma}{p-\sigma}}\int\limits_{w_k < \varepsilon} w_k+\varepsilon^{\frac{p}{p-\sigma}}\int\limits_{w_k \ge \varepsilon} \left(\frac{w_k}{\varepsilon}\right)^{p},
	\end{split}
\end{equation} 
 Choosing a sufficiently small $\varepsilon$, such that $M_2\,\varepsilon^{\frac{\sigma}{p-\sigma}} < \omega_1$,  and subtracting from (\ref{step2}) the first integral of the above right-hand side, multiplied by $M_2$,  the desired inequality (\ref{invPS})  follows with $C=M_2\,\varepsilon^{\frac{p}{p-\sigma}-p}\,\gamma_1^{-1} $.

\vskip0.1cm
In the case  ${\bf A}_{2b}$ ($\max(q,2)-1 \le \sigma \le p-1$ and $\tau>0$) we set $r=\sigma+1$. We can do it by assumption {\bf{(a-e)}} since $\max(q,2)\le r \le p$. This time, being $r-\sigma=1$, the wanted inequality (\ref{invPS})  readily follows from (\ref{step2}) with $C=M_2\gamma_1^{-1}$. We point out that we do not need $\omega_1>0$, but only $\omega_1 \ge 0$, as by default.
\vskip0.1cm
As for  ${\bf A}_{3}$ ($p-1\le  \sigma\le  p$, $\tau >0$ and $\lambda<\infty$), we observe that the case $\sigma=p-1$ has been considered with ${\bf A}_{2b}$, without even assuming $\lambda<\infty$, while the case $\sigma=p$ will be dealt with using a different argument. So we are reduced to $p-1< \sigma < p$, and we turn to set $r=p$ in (\ref{step2}).  Since ${\frac p{p-\sigma}}> p$, using the assumption $\lambda<\infty$, and recalling that $w_k \le \lambda - \lambda_k$,  inequality (\ref{step2}) yields
	\begin{equation}\label{A3}
		\gamma_1\int_\Omega|\nabla w_k|^p \le M_2\int_{\Omega_k} w_k^{s}w_k^{p} \le M_2\lambda^{s}\int_{\Omega_k} w_k^{p},
	\end{equation}
with $s> 0$, so yielding once more inequality (\ref{invPS}), for $C=M_2\lambda^{s}\gamma_1^{-1}$. 

\vskip0.2cm
{\it Case $\sigma \ge p$}
\vskip0.1cm
For $\sigma =p$, which is included in the case ${\bf A}_{3}$ with $\tau >0$ and $\lambda<\infty$, 
we go back to (\ref{step1}) with $r=p$. From this, recalling again that $w_k \le \lambda - \lambda_k$,
	\begin{equation}\label{A3}
	\gamma_1	\int_\Omega|\nabla w_k|^p \le M_2\int_{\Omega} |\nabla w_k|^{p} w_k\le M_2(\lambda-\lambda_k)\int_{\Omega} |\nabla w_k|^p.
	\end{equation}
Here the first integral cannot be zero, otherwise by Poincaré inequality it would be $w_k=0$ a.e., against (\ref{lamda}), as we saw in the case $0<\sigma <p$. It follows that $\gamma_1 \le M_2(\lambda-\lambda_k ) $, which leads to a contradiction for $k$ large, since $\lambda_k \to \lambda$. Therefore $\lambda \le 0$, i.e $u \le v$ in $\Omega$.
\vskip0.2cm 
For $\sigma> p$, recall that we have assumed $u,v \in W^{1,\infty}_{loc}(\Omega)$. Note that, since in particular $u,v \in L^{\infty}_{loc}(\Omega)$, and $u \le v$ on $\partial \Omega$, then $u-v>0$ a.e.  only on sets with positive distance from the boundary, so that  $\lambda < \infty$, and $\lambda-\lambda_k \to 0$ as $k \to \infty$. We still go back to (\ref{step1}), obtaining
	\begin{equation}\label{A3}
		\gamma_0\int_\Omega|\nabla w_k|^p \le M_1\int_{\Omega} |\nabla w_k|^{\sigma}w_k\le (\lambda-\lambda_k) N^{\sigma -p}\int_{\Omega} |\nabla w_k|^p,
	\end{equation}
where we have used the fact that  $w_k=0$ outside a compact subset  $K$ of $\Omega$, taking a constant $N>0$ such that $|\nabla u|+|\nabla v| \le N$ in $K$. As before, the first integral cannot be zero, and so we  are led to a contradiction, $\gamma_0 \le N^{\sigma-p}(\lambda-\lambda_k)$, for $k$ large. We conclude also in this last case that $\lambda \le 0$, i.e. $u \le v$ in $\Omega$, so completing  the proof. \qed
\vskip0.2cm
\begin{rem}{\rm  We remark that, if in the interval $0< \sigma \le 1$, we need in general to assume $\omega>0$ (see ${\bf A}_1$), no additional condition is nevertheless needed when $p=2$ and $\sigma=1$. In fact, looking at the proof of Theorem \ref{cp:thm}: only the basic condition $\tau \ge 0$ needs to have (\ref{step1}) with $r=2$, $\omega_1\ge 0$ and $\sigma=1$, and by Cauchy-Schwarz inequality, we directly obtain (\ref{invPS})  with $r=2$, which is enough, see  (\ref{ctrd}), to prove the comparison principle by contradiction. This is the case of $p$-Laplacian in the linear case, i.e. the Laplacian $\Delta$, and of $(2,q)$-Laplacians $\Delta_{2,q}=\Delta+\Delta_q$, if the Hamiltonian has at most linear growth in the gradient. It continues to hold for $(p,2)$-Laplacians, choosing $r=2$ instead of $r=p$ in the case ${\bf A}_1$.}\qed
\end{rem}
\vskip.0cm
We also give below a proof of Theorem \ref{cp:thmm}, which concerns the case $1 < p < 2$. Here we need $u,v \in W^{1,\infty}_{loc}(\Omega)$ for all $\sigma>0$.
\vskip0.2cm
{\bf Proof of Theorem \ref{cp:thmm}}
\vskip0.1cm
Let $u,v \in W^{1,\infty}_{loc}(\Omega)$ such that $Qu \le Qv$, and $u \le v$ on $\partial \Omega$. As before in the proof of Theorem \ref{cp:thm}, we  argue by contradiction, supposing (\ref{lamda}), i.e. $\lambda=\textup{ess}\sup_\Omega(u-v)>0$. We also keep the same notations, and use the same test functions $w_k$. Let $K$ be a compact subset of $\Omega$ such that  $w_k=0$ outside $K$ for all $k \in \mathbb N$, and $N$ be a positive number such that $|\nabla u|+|\nabla v| \le N$ in $K$. By estimating the left-hand side of (\ref{t1}) with condition  {\bf{(a-e)}}, case $1 < p <2$, and the right-hand side with (\ref{Hxeq}) and the considerations on $\omega$ preceding  (\ref{step0}), we have
	\begin{equation}\label{step00}
		\gamma_1	\int_\Omega|\nabla w_k|^2 +\omega_1\int_\Omega w_k  \le M \int_\Omega ||\nabla u|^\sigma-|\nabla u|^\sigma|w_k
	\end{equation}
with constants $\gamma_1=\gamma_0(1+N)^{p-2}>0$ and $\omega_1 \ge 0$. 
 \vskip0.1cm
 We will show in both cases ${\bf {\widetilde A}}_1$ and ${\bf {\widetilde A}}_w$  an inequality of type (\ref{invPS}) with $r=2$, namely
 	\begin{equation}\label{invPS2}
 		\int_\Omega|\nabla w_k|^2 \le C	\int_{\Omega_k}|w_k|^2,
 	\end{equation}
 with a positive constant $C$ independent of $k$, and we know from the proof of Theorem \ref{cp:thm} that this inequality leads to a contradicion, showing that $\lambda \le 0$, and $u \le v$ in $\Omega$. 
 \vskip0.1cm 
 As in the proof of Theorem \ref{cp:thm}, case $0< \sigma < p$, we are left with proving (\ref{invPS2}).
 \vskip0.1cm
 In the case ${\bf {\widetilde A}}_1$ ($0 < \sigma \le 1$, $\omega>0$), we reason as in Theorem \ref{cp:thm}, case ${\bf {A}}_1$, observing that (\ref{step00}): is nothing else than (\ref{step0}) with $r=2$ (and $\tau_1=0$):  first, we use (\ref{sig1:eq}) to obtain the inequality (\ref{step1}) with $r=2$; next, Young inequality to get (\ref{step2}) with $r=2$; finally, the splitting argument (\ref{sig<1}) with $p=2$ (it is here that condition $\omega>0$ is needed)  to estimate the right-hand side of inequality (\ref{step2}) with $r=2$, and the desired inequality (\ref{invPS2}) follows.
 \vskip0.1cm
 In the case ${\bf {\widetilde A}}_2$ ($ \sigma \ge 1$), we use the inequality
 \begin{equation}\label{1sig}
 	||\xi|^\sigma-|\eta|^\sigma|\le \sigma(|\xi|+|\eta|)^{\sigma-1}|\xi-\eta|,
 \end{equation}
 derived from the  Mean-Value Theorem, to estimate the right-hand side of (\ref{step00}). Here condition $\omega_1\ge0$ is sufficient, and so we get
 	\begin{equation}\label{step01}
 		\gamma_1	\int_\Omega|\nabla w_k|^2  \le \sigma MN^{\sigma-1} \int_\Omega |\nabla w_k| w_k.
 	\end{equation}
 Letting $C=\sigma MN^{\sigma-1}\gamma_1^{-1}$, by using  Cauchy-Schwarz inequality in the right-hand side, we get
 	\begin{equation}\label{preinvPS2}
 \int_\Omega|\nabla w_k|^2  \le C  \left(\int_\Omega|\nabla w_k|^2\right) ^{1/2} \left(\int_\Omega w_k^2\right) ^{1/2}.
 \end{equation}
 From this, since we know the first integral has to be positive (see the proof of Theorem \ref{cp:thm}, case $\sigma \ge p$), we get (\ref{invPS2}), and this concludes the proof.\qed
  \vskip0.3cm 
  {\bf Outline of the proofs of Theorems \ref{cpfg:thm} and \ref{cpfg:thmm}} 
	\vskip0.1cm
The proof of Theorems \ref{cpfg:thm} and \ref{cpfg:thmm} follows the same lines of the respective proofs of Theorems \ref{cp:thm} and \ref{cp:thmm}. Comparing (\ref{fgxeq}), (\ref{fpos}) and (\ref{gpos}) with (\ref{Hxeq}), (\ref{omegapos}) and (\ref{taupos}), respectively, in the previous Section, it is sufficient to replace $\omega(x,u-v)$ and $\tau(x,u-v)$ with $f(x,u)-f(x,v)$ and $g(x,u)-g(x,v)$, respectively. Likewise, taking into account that $u,v \in L^\infty_{loc}(\Omega)$ in both Theorems  \ref{cp:thm} and \ref{cp:thmm}, we take $M=\max_{K \times J}|g|$. Here $K$ is the compact subset of $\Omega$ used in the proof of Theorem \ref{cp:thm}, and  $J=[a,b]$ with $a=\min(\text {ess} \,\inf_Ku,\text{ess}\, \inf_Kv)$ and $b=\max(\text{ess}\, \sup_Ku,\text{ess}\, \sup_Kv)$. The remainder of the proof flows verbatim.\qed

\begin{rem} 
	Condition $\tau>0$, assumed in the cases ${\bf A}_{2a}$, ${\bf A}_{2b}$, ${\bf A}_{3}$  and ${\bf A}_{4}$ of Theorem \ref{cp:thm}, can be replaced by condition $\omega>0$, provided that either $u$ or $v$ belong to $W^{1,\infty}_{loc}(\Omega)$. In fact, suppose that $u\in W^{1,\infty}_{loc}(\Omega)$, and $|\nabla u|\le N$ in the compact $K\subset \Omega$ which contains the supports of all $w_k$: here we have to estimate the second integral of the right-hand side in (\ref{tauest}) choosing $\varepsilon >0$ such that $\frac M2\,((1+\varepsilon)^{\sigma-1}-1)N^\sigma< \omega_1$. In this way, subtracting this integral from (\ref{tauest}), we obtain (\ref{step1}) with a positive constant $\omega_2<\omega_1$, and the remainder of the proof can be carried out with $\omega_2$ instead of $\omega_1$. By analogy, the condition that $g$ is increasing, assumed in the cases ${\bf A}_{2a}$, ${\bf A}_{2b}$, ${\bf A}_{3}$ and ${\bf A}_{4}$ of Theorem \ref{cpfg:thm}, can be replaced by the condition that $f$ is increasing.\qed
\end{rem}

\begin{rem} 
	Alternatively, in the case ${\bf A}_{3}$ of Theorem \ref{cpfg:thm}, when $p-1<\sigma \le p$, we can drop the condition that $g$ is increasing if we assume $g(x,v)\ge0$. This result is known in the case of $p$-Laplacian, see \cite[Theorem 1.5]{MOV}. In Theorem \ref{g>0} below, we extend it to our more general setting based on condition {\bf{(a-e)}}, but unfortunately, requiring an homogeneity assumption on the vector field $\bm a$, our result does not hold for the $(p,q)$-Laplacian. We also provide a variation of the proof , which allows to improve the result even in the case of $p$-Laplacian, assuming $u^-\in L^\infty(\Omega)$ rather than $u \in L^\infty(\Omega)$.
	\qed
\end{rem}

\begin{thm}\label{g>0}
	Assume $2 \le  p<\infty$, $q=p$ and $p-1<\sigma\le p$. Let condition {\bf{(a-e)}} be satisfied, and $f(x,t)$, $g(x,t)$ be non-decreasing function in $t \in \R$ for all $x \in \Omega$. Suppose  $u,v \in W^{1,p}_{loc}(\Omega)\cap L^\infty_{loc}(\Omega)$, $u$ essentially bounded from below and  $g(x,v)$ non-negative in $\Omega$. In addition, suppose that the following homogeneity conditions hold:
	\vskip0.2cm
	\begin{itemize}
		\item[{\bf (a-h)}] \quad 	$\bm{a}(x,\delta\xi) =\delta^{p-1} \bm{a}(x,\xi)$, \;\; $(x,\xi) \in \Omega \times \R^n$

	\vskip0.2cm
		\item[{\bf (f-h)}] \quad 	$f(x,\delta t) \le \delta^{p-1} f(x,t)$, \;\; $(x,t) \in \Omega \times \R^n$
			\end{itemize}
			\vskip0.2cm
			
\noindent 			for all $\delta\in(0,1)$. If $Q^*u \le 0 \le Q^*v$ in $\Omega$, and $u \le v$ on $\partial\Omega$, then $u \le v$ in $\Omega$
	 
\end{thm}
{\bf Proof.} As before,  in the proof of Theorem \ref{cp:thm}, we argue by contradiction, supposing on the contrary (\ref{lamda}), namely $\lambda=\text{ess} \sup_\Omega (u-v)>0$. Recall that by condition $u \le v$ on $\partial \Omega$, since $u,v \in L^\infty_{loc}(\Omega)$, we have $\lambda< \infty$. 

Next, we introduce a further parameter $\varepsilon >0$, in order to penalize the function $u$, bounded below, setting 
\begin{equation}\label{ueps}
	u_\varepsilon(x)= u(x) -\varepsilon(u(x)+\underline \mu),
\end{equation}
where $\underline \mu$ is a positive constant such that $\text{ess}\inf_\Omega u \ge - \underline \mu$. Correspondingly, we define
\begin{equation}\label{lamdaeps}
	\lambda_{\varepsilon} = \text{ess} \sup_\Omega (u_\varepsilon -v), 
\end{equation}
and we note that, being $u_\varepsilon < u$, that $\lambda_\varepsilon \le \lambda < \infty$. Taking $\varepsilon>0$ eventually smaller, we can also keep $\lambda_\varepsilon>0$. In fact, by condition $u \le v$ on $\partial \Omega$, for a positive number $\ell<\lambda$, the set $O_\ell=\{x \in \Omega: u(x)-v(x)> \ell\}$ has a positive distance from $\partial \Omega$, so that there is a compact $K_\ell\subset \Omega$ such that $O_\ell\subset K_\ell$. Moreover, by definition of $\lambda$, the set $O_\ell$ has positive measure, and in $O_\ell$, recalling that $u \in L^\infty_{loc}(\Omega)$, we have 
$$u_\varepsilon-v = u-\varepsilon(u+ \underline\mu)-v > \ell-\varepsilon(\overline \mu+\underline\mu), $$
where $\overline \mu$ is a positive number greater than  $\text{ess}\sup_{K_\ell}u$. Therefore, choosing $\varepsilon>0$ such that $\ell> \varepsilon (\overline \mu+\underline\mu)$, 
we have $\text{ess} \sup_{O_\ell}(u_\varepsilon -v )>0$, and this shows that $\lambda_\varepsilon >0$ as claimed.

Then we can take an increasing sequence of positive numbers $\lambda_{\varepsilon,k}$ such that $\lambda_{\varepsilon,k} \to \lambda_\varepsilon$ as $k \to \infty$. Let $\mathcal O_{\varepsilon,k}=\{x \in \Omega : u_\varepsilon-v>\lambda_{\varepsilon,k}\} $. If $x \in \mathcal O_{\varepsilon,k}$, then $u(x)-v(x)> \lambda_{\varepsilon,k} > \lambda_{\varepsilon,1}$, and so, by virtue of the boundary condition $u \le v$ on $\partial \Omega$, all $\mathcal O_{\varepsilon,k}$ can be included in a compact subset $K\subset \Omega$.  Consequently, the functions $w_{\varepsilon,k}:=(u_\varepsilon-v-\lambda_{k,\varepsilon})^+\in W^{1,p}(\Omega)$ have support in $K\subset \Omega$, and so they can be used as test functions for the inequality $Q^*u \le 0$, namely:
	\begin{equation}\label{t1u}
	\int_\Omega \bm{a}(x,\nabla u)\cdot \nabla w_{\varepsilon,k}\le 	-\int_\Omega(f(x,u)+g(x,u)|\nabla u|^\sigma)w_{\varepsilon,k}.
\end{equation}
Taking into account that $u=\frac{u_\varepsilon +\varepsilon \underline \mu}{1-\varepsilon}$, by assumptions on $\bm a$, $f$ and $g$, we get 
	\begin{equation}\label{t2u}
	\int_\Omega \bm{a}(x,\nabla u_\varepsilon)\cdot \nabla w_{\varepsilon,k}\le-\int_\Omega(1-\varepsilon)^{p-1}\left(f(x,\tfrac{u_\varepsilon}{1-\varepsilon})+(1-\varepsilon)^{-\sigma}g(x,\tfrac{u_\varepsilon}{1-\varepsilon}\right)|\nabla u_\varepsilon|^\sigma)w_{\varepsilon,k}.
\end{equation}
Using the same functions $w_{\varepsilon,k}$ to test the inequality $Q^*v \ge 0$, we get
	\begin{equation}\label{t1v}
	-\int_\Omega \bm{a}(x,\nabla v)\cdot \nabla w_{\varepsilon,k}\le 	\int_\Omega(f(x,v)+g(x,v)|\nabla v|^\sigma)w_{\varepsilon,k}.
\end{equation}
Summing (\ref{t2u}) and (\ref{t1v}), we have therefore 
	\begin{equation}\label{t3}
	\int_\Omega (\bm{a}(x,\nabla u_\varepsilon)-\bm{a}(x,\nabla v))\cdot \nabla w_{\varepsilon,k}
	\le \int_\Omega \left(g(x,v)|\nabla v|^\sigma-(1-\varepsilon)^{-\sigma+p-1}g(x,\tfrac{u_\varepsilon}{1-\varepsilon})|\nabla u_\varepsilon|^\sigma\right)w_{\varepsilon,k}
\end{equation}
where we have used the fact that, by {\bf{(f-h)}} the non-decreasing monotonicity of $f$,
\begin{equation}\label{finc}
	f(x,v)-(1-\varepsilon)^{p-1}f(x,\tfrac{u_\varepsilon}{1-\varepsilon})\le f(x,v)-f(x,u_\varepsilon) \le0 
\end{equation}
on the support of $w_{\varepsilon,k}$, where $u_\varepsilon-v\ge \lambda_{\varepsilon,k}$.

Now, we estimate from below the left-hand side of (\ref{t3}) by condition {\bf{(a-e)}} with $r=p=q$, obtaining
	\begin{equation}\label{t4}
		\gamma_0\int_\Omega |\nabla w_{\varepsilon,k}|^p\le 
		\int_\Omega \left(g(x,v)|\nabla v|^\sigma-(1-\varepsilon)^{-\sigma+p-1}g(x,\tfrac{u_\varepsilon}{1-\varepsilon})|\nabla u_\varepsilon|^\sigma\right)w_{\varepsilon,k}.
\end{equation}
Next, we majorize the right-hand side. Here we use the non-decreasing monotonicity of $g$, which yields, still observing that $v \le u_\varepsilon$ on the support of $w_{\varepsilon,k}$,
\begin{equation}\label{ginc}
	g(x,v)|\nabla v|^\sigma-(1-\varepsilon)^{-\sigma+p-1}g(x,\tfrac{u_\varepsilon}{1-\varepsilon})|\nabla u_\varepsilon|^\sigma\le g(x,v)\left(|\nabla v|^\sigma-\frac{|\nabla u_\varepsilon|^\sigma}{(1-\varepsilon)^{\sigma-p+1}}\right),
\end{equation} 
and the inequality 
\begin{equation}\label{A4}
	|\xi|^\sigma-\frac{|\eta|^\sigma}{(1-\varepsilon)^{\sigma-p+1}} \le C_\varepsilon|\xi-\eta|^\sigma, \; \; \xi,\eta \in \R^n,
\end{equation}
where $C_\varepsilon$ is a positive constant (see \cite[Lemma A.4]{MOV} for a proof), to get
	\begin{equation}\label{t5}
	\gamma_0\int_\Omega |\nabla w_{\varepsilon,k}|^p\le M_1^* \int_\Omega |\nabla w_{k,\varepsilon}|^\sigma w_{\varepsilon,k},
\end{equation}
where $M_1^*=C_\varepsilon \max_{_{K\times J_v} }g$, where $J_v=[\text{ess}\inf_K v, \text{ess}\sup_K v]$.

For $\sigma \in (p-1,p)$, using once more  Young inequality as in (\ref{young}), we obtain
	\begin{equation}\label{t6}
	\gamma_1\int_\Omega |\nabla w_{\varepsilon,k}|^p\le M_2^*\int_{\Omega} w_{\varepsilon,k}^{\frac p{p-\sigma}},
\end{equation}
with positive constants $\gamma_1$ and $M_2^*$, independent of $k$. From this, since $w_{\varepsilon,k} \le \lambda_\varepsilon -\lambda_{\varepsilon,k}$, and $w_{k,\varepsilon}^{\frac p{p-\sigma}} =w_{k,\varepsilon}^sw_{k,\varepsilon}^p$, where $s=\frac{p}{p-\sigma}-p$ is a positive number, we deduce that 
\begin{equation}\label{invP}
	\int_\Omega  |\nabla w_{\varepsilon,k}|^p\le C(\lambda_\varepsilon-\lambda_{\varepsilon,k})^s\int_{\Omega} w_{\varepsilon,k}^p
\end{equation}
for a positive constant $C$. Combined with the Poincaré inequality, this yields 
\begin{equation}\label{Poinc}
	\int_\Omega w_{\varepsilon,k}^p \le C\int_\Omega|\nabla w_{\varepsilon,k}|^p  \le C(\lambda-\lambda_k)^s \int_{\Omega} w_{\varepsilon,k}^p,
\end{equation}
where $C$ is a positive constant with possibly different values, independent of $k$. In the above we cannot have $w_{\varepsilon,k}=0$ a.e. in $\Omega$ for any $k \in \mathbb N$, otherwise $u_\varepsilon-v \le \lambda_{\varepsilon,k}$ a.e. in $\Omega$, which would contradicts (\ref{lamdaeps}). Then we get the inequality $1\le C(\lambda_{\varepsilon}-\lambda_{\varepsilon,k})$, which leads to a contradiction for $k$ large, since $\lambda_{\varepsilon,k}\to \lambda_{\varepsilon}$ as $k \to \infty$. We conclude therefore that  $\lambda \le 0$, namely $u \le v$ in $\Omega$.

If $\sigma=p$, then if follows directly from (\ref{t5}) that 
	\begin{equation}\label{t5'}
	\gamma_0\int_\Omega |\nabla w_{\varepsilon,k}|^p\le M_1^* (\lambda_\varepsilon-\lambda_{\varepsilon_k})\int_\Omega |\nabla w_{k,\varepsilon}|^p w_{\varepsilon,k},
\end{equation}
where  the integral cannot be zero otherwise, by Poincaré inequality, $w_{\varepsilon,k}=0$ a.e., which would contradicts (\ref{lamdaeps}), as seen in the previous case. We deduce that $\gamma_0 \le M_1^*(\lambda_\varepsilon-\lambda_{\varepsilon_k})$, which leads to a contradiction for $k$ large. As before,  we conclude that $\lambda \le 0$, and so $u \le v$ in $\Omega$, finishing the proof.
\qed

\end{document}